\documentclass{article}
\usepackage{amsmath}
\usepackage{amssymb}
\usepackage{amsthm}
\usepackage{amsxtra}
\usepackage{pstcol}
\usepackage{graphicx}

\begin{document}
\centerline{\bf Counting numerical semigroups}
\vskip5mm
\centerline{E. Kunz and R. Waldi}
\vskip 3mm\noindent
{\bf Abstract.} We are interested in formulas for the number of elements in certain classes of numerical semigroups
\vskip3mm\noindent
\vskip3mm\noindent
{\bf 1. Introduction}
\vskip3mm

For an integer $p\ge 3$ let $\mathfrak H_p$ be the set of all numerical semigroups containing $p$. Using the notion of Ap\'ery set one can construct a bijective map from $\mathfrak H_p$ onto the set of all lattice points of a certain polyhedral cone $C_p\subset \mathbb R^{p-1}$ of dimension $p-1$ ([Ku],[RGGB]). The lattice points in the interior $C_p^0$ of $C_p$ are in one-to-one correspondence with the $H\in\mathfrak H_p$ of maximal embedding dimension $p$. Ap\'ery's characterization of symmetric semigroups [A] allows to show that the symmetric $H\in\mathfrak H_p$ are mapped onto the lattice points of certain $\lfloor\frac{p}{2}\rfloor$-dimensional closed faces of $C_p$. Here we describe also the distribution of the lattice points corresponding to the pseudo-symmetric semigroups (Proposition 3.2). Further semigroups correspond to the intersection of $C_p$ or its faces with hyperplanes, hence with the lattice points in polyhedrons, and their number can be expressed by quasi-polynomials (Theorem 4.2).

This is the case for the $H\in\mathfrak H_p$ with fixed genus $g$. Also the semigroups $H\in\mathfrak H_p$ containing another number $q$ which is prime to $p$ can be described using hyperplane sections of $C_p$ (Example 4.1b)). We denote this set of numerical semigroups by $\mathfrak H_{pq}$.

We want to study the degree, the leading term and a quasi-period of the involved Ehrhart quasi-polynomials. For the semigroups $H\in \mathfrak H_{pq}$ the leading term is constant and gives therefore an asymptotic estimate for $q\to \infty$ of the number of the $H\in \mathfrak H_{pq}$ (Proposition 5.2). Similarly for the $H\in \mathfrak H_p$ of maximal embedding dimension $p$ and the symmetric $H\in \mathfrak H_{pq}$ (Propositions 4.5 and 4.6).
\hfill\eject
The following figure illustrates the situation in the simplest case $p=3$.

\pspicture*(-1.1,-1.1)(7.1,10.6)
\psgrid[gridlabels=0,subgriddiv=2,gridwidth=0.001pt,subgridwidth=0.001pt,gridcolor=lightgray](-1.1,-1.1)(7.1,10.6)
\psline[linewidth=.5pt](0,0)(0,10.3)
\psline[linewidth=.5pt](0,0)(7,0)
\psdots*[dotsize=3pt](0,0)(1,0)
\psdots*[dotsize=3pt](1,1)(2,1)(3,1)
\psdots*[dotsize=3pt](1,2)(2,2)(3,2)(4,2)(5,2)
\psdots*[dotsize=3pt](2,3)(3,3)(4,3)(5,3)
\psdots*[dotsize=3pt](2,4)(3,4)(4,4)(5,4)
\psdots*[dotsize=3pt](3,5)(4,5)(5,5)
\psdots*[dotsize=3pt](3,6)(4,6)(5,6)
\psdots*[dotsize=3pt](4,7)(5,7)
\psdots*[dotsize=3pt](4,8)(5,8)
\psdots*[dotsize=3pt](5,9)
\psdots*[dotsize=3pt](5,10)
\psline[linewidth=.5pt](-0.33,-0.67)(5.25,10.5)
\psline[linewidth=.5pt](-0.33,-0.67)(6,2.5)
\psline[linestyle=dashed,
linewidth=.5pt](.2,.1)(6,3)
\psline[linestyle=dashed, linewidth=.5pt](.5,0)(5.25,9.5)
\psline[linestyle=dotted, linewidth=1.2pt](2.4,4.6)(5,2)
\uput{3pt}[180](-.4,-.6){$v$}
\uput{2pt}[0](5.2,5){$C_3$}
\uput{2pt}[180](2,1.1){$\delta_1$}
\uput{2pt}[180](1,2){$\delta_2$}

\endpspicture
\vskip1mm

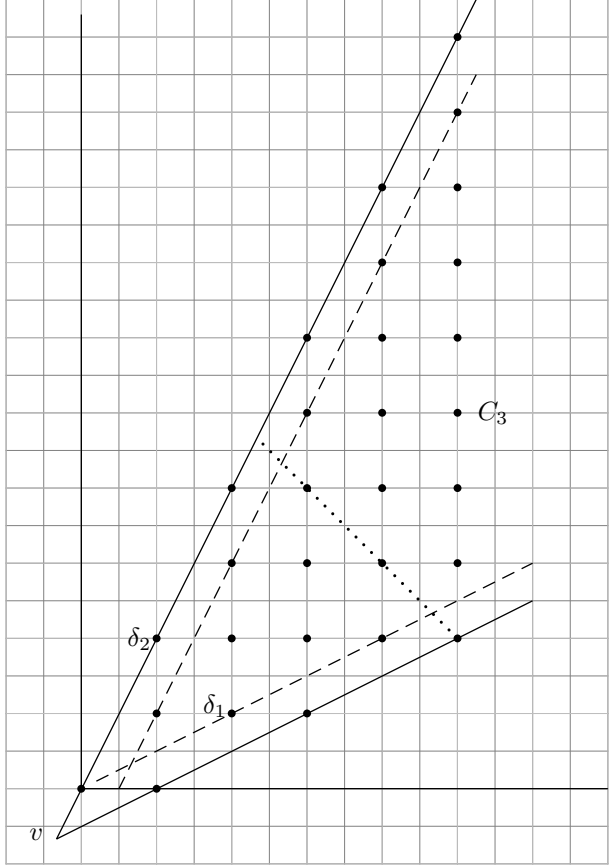
Figure 1: Geographical distribution of the $H\in \mathfrak H_3$.
\vskip2mm

{\footnotesize{The symmetric semigroups (here=complete intersections) belong to the edges of $C_3$ (see Proposition 3.1), those with maximal embedding dimension 3 (here=almost complete intersections) to the interior of $C_3$. The living space of the pseudo-symmetric semigroups are the dashed lines $X_1=2X_2$ and $X_2=2X_1-1$ (see Proposition 3.2) where the origin (which corresponds to $H=\mathbb N$) has to be excluded. That of the semigroups with genus $g=7$ is the dotted line $X_1+X_2=7$. Parallels to this line correspond to semigroups with fixed genus. The figure shows, for example, that for given $g\in\mathbb N$ there is exactly one symmetric (complete intersection) semigroup in $\mathfrak H_3$ of genus $g$ if $g\not\equiv 2\ \text{mod}3$, and no such semigroup if $g\equiv 2\ \text{mod}3$. Similarly there is exactly one pseudo-symmetric $H\in \mathfrak H_3$ of genus $g>0$ if $g\not\equiv 1\ \text{mod}3$, and no such semigroup if $g\equiv 1\ \text{mod} 3$.
The lattice points of $C_3$ with $X$-coordinate $\le i$ (with $Y$-coordinate $\le i$) correspond to the $H\in \mathfrak H_3$ with $1+3i\in H$ (with $2+3i\in H$), see Example 4.1b).}}

Using another visualization of the $H\in \mathfrak H_{pq}$ by lattice paths in the plane  ([KKW],[KW]) recursion formulas for the number of semigroups $H\in \mathfrak H_{pq}$ which interest us can be derived, see Section 5. Finally in Section 6 explicit formulas are given for $p\le 5$. For related counting problems, see [Ka], [BGP].
\vskip3mm
\noindent
{\bf 2. The polyhedral cone $C_p$ and its faces}
\vskip3mm
We recall some facts about $C_p$ which are relevant for us. Details can be found in [Ku] or [RGGB]. For $H\in \mathfrak H_p$ let
$\{h_1,\dots,h_{p-1}\}=\text{Ap}(H,p)$ be the Ap\'ery set of $H$ with respect to
$p$, that is, $h_i$ is the smallest element of $H$ in $i+p\mathbb N$
for $i=1,\dots,p-1$. Write $h_i=i+\mu_ip$. The semigroup $H$ is uniquely
determined by $(\mu_1,\dots,\mu_{p-1})$ since
$H=<p,h_1,\dots,h_{p-1}>$. The points $(\mu_1,\dots,\mu_{p-1})\in
\mathbb N^{p-1}$ are the solutions in $\mathbb N^{p-1}$ of the
system of linear inequalities
$$\begin{cases} X_i+X_j\ge X_{i+j}\qquad\qquad (i+j<p)\cr
X_i+X_j\ge X_{i+j-p}-1\quad\ (i+j>p)\end{cases}\leqno(1)$$
Let $C_p$ be the solution set of (1) in $\mathbb R^{p-1}$.
This is a polyhedral cone with vertex $v:=(-\frac{1}{p},-\frac{2}{p}\dots,-\frac{p-1}{p})$ and
$$\mu:\mathfrak H_p\to C_p\cap\mathbb N^{p-1}\quad (H\mapsto (\mu_1,\dots,\mu_{p-1}))$$
is a bijection of $\mathfrak H_p$ onto the set of lattice points of $C_p$. For $H\in \mathfrak H_p$ we call $\mu(H)$ the lattice point associated to $H$, and for $P=\mu(H)$ we say that $H$ is the semigroup belonging to the point $P$.
By a {\it face} of $C_p$ we always understand an open face.
Its dimension is the dimension of the smallest affine space containing it.
We consider the interior $C_p^0$ of $C_p$ and the vertex $v$ also as faces.
\vskip2mm
We have $C_p^0\cap\mathbb N^{p-1}=(1,\dots,1)+C_p\cap \mathbb N^{p-1}$ ([Ku], 1.4d).
$C_p$ and $C_p^0$ have dimension $p-1$. The {\it facets} (1-codimensional faces) correspond bijectively to the hyperplanes
$$\begin{cases} E_{ij}: X_i+X_j-X_{i+j}=0\qquad\qquad (i+j<p)\cr
E_{ij}: X_i+X_j-X_{i+j-p}+1=0\quad\ (i+j>p)\end{cases}$$
Hence for even $p$ there are $\frac{(p-1)^2-1}{2}$ facets, for odd $p$ their number is $\frac{(p-1)^2}{2}$.

The translation of $C_p$ by the vector $-v$ leads to a polyhedral cone $C_p^*:=-v+C_p$ in $\mathbb R_+^{p-1}$ with vertex at the origin. It is the solution set of the system $X_i+X_j\ge X_{i+j}\ (i+j\ne p\ \text{with indices reduced modulo}\ p)$. Therefore $C_p^*\subset C_p$. On each face of $C^*_p$ there are lattice points. This is true in particular for its edges (1-dimensional faces). On each such edge there is a {\it primitive lattice point} $\delta$, that is a point with relatively prime integral coordinates, and all other lattice points on the edge are integral multiples of $\delta$. The set of all these $\delta$ is called the {\it canonical system of representatives} for the edges of $C_p^*$.

Let $m(H)$ denote the multiplicity and edim$(H)$ the embedding dimension of a numerical semigroup $H$. The lattice points in the interior of $C_p$ correspond to the $H\in \mathfrak H_p$ with edim($H)=m(H)=p$ ([Ku], 2.4b), called in modern language the {\it semigroups of maximal embedding dimension} $p$.

Let $S$ be a $d$-dimensional face of $C_p$, $S^*$ the face of $C_p^*$ parallel to $S$ and $\bar S$ resp. $\bar S^*$ the topological closures of these faces. The generating function (Hilbert series) for the lattice points on $S$ is the formal power series
$$H_S(X_1,\dots,X_{p-1})=\sum_{\mu\in S\cap\mathbb N^{p-1}}X^{\mu}.$$
The generating functions $H_{\bar S}$ for $\bar S$, $H_p$ for $C_p$ and $H_p^0$ for the interior $C_p^0$ of $C_p$ are defined accordingly.
\vskip3mm\noindent
{\bf 2.1. Remark.} ([Ku], 2.1) $H^0_p(X_1,\dots,X_{p-1})=X_1\cdots X_{p-1}H_p(X_1,\dots,X_{p-1})$.
\vskip3mm

Let $\{\delta_1,\dots,\delta_t\}$ be the
canonical system of representatives for the edges of $C_p^*$ and
$\delta_1,\dots,\delta_r$ the vectors belonging to $\bar{S^*}$. Using slack variables $T_{ij}$ the system (1) of Section 1 becomes a linear system of equations
$$X_i+X_j-X_{i+j}-T_{ij}=0\ (i+j<p)$$
$$X_i+X_j-X_{i+j-p}-T_{ij}=-1\ (i+j>p).$$
To it the results of Stanley [S1] Chap.I can be applied. One obtains
\vskip2mm \noindent
{\bf 2.2. Proposition.} $H_{\bar S}$ can be written in the form
$$H_{\bar S}=\frac{F_{\bar S}}{\prod_{i=1}^r(1-X^{\delta_i})},\ F_{\bar S}\in \mathbb Z[X_1,\dots,X_{p-1}]$$
Similarly $H_p$ has the form
$$H_p=\frac{F_p}{\prod_{i=1}^t(1-X^{\delta_i})},\ F_p\in \mathbb Z[X_1,\dots,X_{p-1}].$$
Further deg$_{X_i}H_{\bar S}\le -1$ and deg$_{X_i}H_p\le -1\ (i=1,\dots,p-1)$.
\vskip2mm

See [Ku], 3.2 and 3.3 for details. The proof of [Ku], 3.2 contains the false statement that $R_{\overline{S^*}}:=\mathbb C[\{X^{\mu}\}_{\mu\in \overline{S^*}\cap\mathbb N^{p-1}}]$ as an algebra over $\mathbb C$ is generated by $X^{\delta_1},\dots,X^{\delta_r}$. However the algebra is a finitely generated module over $A:=\mathbb C[X^{\delta_1},\dots,X^{\delta_r}]$. Then also the module $M_{\bar S}:=\bigoplus_{\mu\in \bar S\cap \mathbb N^{p-1}}\mathbb C\cdot X^{\mu}$ is a finitely generated $A$-module, which is what is actually used in the proof.
\vskip3mm\noindent
{\bf 3. Symmetric and pseudo-symmetric semigroups}
\vskip3mm
Ap\'ery's characterization of symmetric semigroups in terms of their Ap\'ery sets ([A]) allows to describe the lattice points in $C_p$ belonging to such semigroups. Let
$$t:=\begin{cases}p-1,\ p\ \text{odd}\cr \frac{p}{2},\qquad p\ \text{even}\end{cases}.$$
\noindent
{\bf 3.1. Proposition} ([Ku], 2.9) There are $t$ faces $S_1,\dots,S_t$ of dimension $\lfloor\frac{p}{2}\rfloor$ such that any symmetric $H\in\mathfrak H_p$ belongs to exactly one $\bar S_j$, and all lattice points on the $\bar S_j$ correspond to symmetric $H\in \mathfrak H_p$.
\vskip2mm
In particular $\bar S_j\cap \bar S_k\cap \mathbb N^{p-1}=\emptyset$ for $j\ne k$. According to Proposition 2.2 the generating function for the symmetric semigroups $$H_{\text{sym}}(X_1,\dots,X_{p-1}):=\sum_{\mu\in \cup_{j=1}^t\bar S_j\cap\mathbb N^{p-1}}X^{\mu}=\sum_{j=1}^tH_{\bar S_j}(X_1,\dots,X_{p-1})$$
is a rational function.

For pseudo-symmetric semigroups the situation is more complicated. Remember that a numerical semigroup $H$ of genus $g(H)$ and Frobenius number $F(H)$ is {\it pseudo-symmetric} if $2g(H)=F(H)+2$.

Let $\mathfrak S_{p-1}$ be the permutation group of $\{1,\dots,p-1\}$ and
$\mathfrak S_{p-1}^*$ the set of all $\sigma \in \mathfrak S_{p-1}$ such that
$$\sigma(i)+\sigma(p-2-i)\equiv \sigma(p-2)\ \text{mod}p\ (i=1,\dots,p-3)$$
$$2\sigma(p-1)\equiv\sigma(p-2)\ \text{mod}p.$$
For $\sigma\in \mathfrak S_{p-1}^*$ let $L_{\sigma}$ be the affine subspace of $\mathbb R^{p-1}$ defined by the linear equations
$$X_{\sigma(i)}+X_{\sigma(p-2-i)}=\begin{cases}X_{\sigma(p-2)}\qquad \text{if}\ \sigma(i)+\sigma(p-2-i)=\sigma(p-2)\cr
X_{\sigma(p-2)}-1\ \text{if}\ \sigma(i)+\sigma(p-2-i)=\sigma(p-2)+p\end{cases}\leqno (1)$$
$(i=1,\dots,p-3)$ and
$$2X_{\sigma(p-1)}=X_{\sigma(p-2)}+\begin{cases}1\quad \text{if}\ 2\sigma(p-1)=\sigma(p-2)\cr
0\quad \text{if}\ 2\sigma(p-1)=\sigma(p-2)+p\end{cases}.$$
This space has dimension $\lfloor\frac{p-1}{2}\rfloor$.

The following Proposition is based on [RG],4.15.
\vskip3mm\noindent
{\bf 3.2. Proposition.} a) For $\sigma\in \mathfrak S_{p-1}^*$ each lattice point $(\mu_{\sigma(1)},\dots,\mu_{\sigma(p-1)})\ne 0$ of $L_{\sigma}\cap C_p$ belongs to a pseudo-symmetric semigroup.

\noindent
b) For each pseudo-symmetric $H\in \mathfrak H_p$ there exists a
$\sigma\in \mathfrak S_{p-1}^*$ such that $\mu(H)\in L_{\sigma}\cap C_p$.

\noindent
c) For $p>3$ the $L_{\sigma}\cap C_p$ are contained in the boundary of $C_p$.

\noindent
d) $L_{id}\cap C_p=\bar S\cap H$ where $\bar S$ is the closure of a face $S$ of $C_p$ of dimension $\lfloor\frac{p+1}{2}\rfloor$ containing $0$ and $H$ a hyperplane through $0$ defining a facet of $C_p^*$.
\vskip3mm\noindent
{\bf Proof.} a) Let $h_{\sigma(i)}:=\sigma(i)+\mu_{\sigma(i)}p\ (i=1,\dots,p-1)$. Then $\{h_{\sigma(1)},\dots,h_{\sigma(p-1)}\}$ is the Ap\'ery set of a semigroup $H\in \mathfrak H_p$ and $\sum_{i=1}^{p-1}h_{\sigma(i)}={{p}\choose{2}}+p\sum_{i=1}^{p-1}\mu_{\sigma(i)}={{p}\choose{2}}+g(H)p$.

From the equations defining $L_{\sigma}$ we obtain
$$h_{\sigma(i)}+h_{\sigma(p-2-i)}=h_{\sigma(p-2)}\ (i=1,\dots,p-3)$$
$$2h_{\sigma(p-1)}=h_{\sigma(p-2)}+p.$$
Obviously $h_{\sigma(i)}<h_{\sigma(p-2)}$ for $i=1,\dots,p-3$, hence $h_{\sigma(p-2)}\ge p-2$, and then by the last equation
$h_{\sigma(p-1)}\ge p-1$ and $h_{\sigma(p-2)}-h_{\sigma(p-1)}=h_{\sigma(p-1)}-p\ge -1$. In case $h_{\sigma(p-1)}>h_{\sigma(p-2)}$ we would have $h_{\sigma(p-1)}=p-1, \mu_{\sigma(p-1)}=0$. But then $\mu_i=0$ for $i=1,\dots,p-1$, contrary to the assumption.

Thus $h_{\sigma(p-2)}$ is the maximal element of the Ap\'ery set of $H$, i.e. $h_{\sigma(p-2)}=F(H)+p$ with the Frobenius number $F(H)$ of $H$.
Adding the above equations gives
$$2\sum_{i=1}^{p-1}h_{\sigma(i)}=ph_{\sigma(p-2})+p=p(F(H)+p)+p$$
that is
$$2g(H)=F(H)+2$$
which is equivalent to $H$ being pseudo-symmetric.

\noindent
b) If $H$ is pseudo-symmetric, then $F(H)$ is an even number. By [RG], 4.15 there is a permutation $\sigma\in \mathfrak S_{p-1}^*$ such that $\{h_{\sigma(1)},\dots,h_{\sigma(p-1)}\}$ is the Ap\'ery set of $H$ and
$$h_{\sigma(i)}+h_{\sigma(p-2-i)}=h_{\sigma(p-2)}\ (i=1,\dots,p-3)$$
$$h_{\sigma(p-2)}=F(H)+p, h_{\sigma(p-1)}=\frac{F(H)}{2}+p$$
in particular
$$2h_{\sigma(p-1)}=F(H)+2p=h_{\sigma(p-2)}+p.$$
With $h_{\sigma(i)}:=\sigma(i)+\mu_{\sigma(i)}p\ (i=1,\dots,p-1)$ it follows that $(\mu_{\sigma(1)},\dots,\mu_{\sigma(p-1)})\in L_{\sigma}\cap C_p$.

\noindent
c) For $p>3$ the equations (1) show that $L_{\sigma}$ is contained in a hyperplane which defines a facet of $C_p$, hence $L_{\sigma}\cap C_p$ belongs to the boundary of $C_p$.

\noindent
d) $\sigma=id$ belongs to $\mathfrak S_{p-1}^*$ and in this case $L_{\sigma}$ is the intersection of the hyperplanes
$$H_i: X_i+X_{p-2-i}=X_{p-2}\ (i=1,\dots,p-3), H: 2X_{p-1}=X_{p-2}.$$
The $H_i$ define facets of $C_p$ containing the origin, and $H$ a facet of $C_p^*$. We have
$$L_{id}\cap C_p=\cap_{i=1}^{p-3}(H_i\cap C_p)\cap H$$
and $\bar S:=\cap_{i=1}^{p-3}(H_i\cap C_p)$ is the closure of a face $S$ of $C_p$.

The point $P:=(1,\dots,1,2,1)$ is in $L_{id}\cap C_p$ and in the interior of the half-spaces $X_i+X_j\ge X_{i+j} (i+j\ne p-2)$. Therefore an open neighborhood of $P$ in the $\lfloor\frac{p-1}{2}\rfloor$-dimensional
affine space $L_{id}$ is contained in $L_{id}\cap C_p$. Hence the face $S$ must have dimension $\lfloor\frac{p+1}{2}\rfloor$.\ $\Box$
\vskip3mm\noindent
{\bf 3.3. Corollary.} ([RG], 4.26). For $p>3$ pseudo-symmetric semigroups of $\mathfrak H_p$ have embedding dimension $<p$.
\vskip3mm\noindent
{\bf Proof.} Otherwise the lattice points of such semigroups would be in the interior of $C_p$ contradicting 3.2c).\ $\Box$
\vskip3mm
\noindent
{\bf 3.4. Example.} In general $L_{id}\cap C_p$ is not a polyhedral cone. If $p=7$, then $L_{id}$ is a 3-space, in which $L_{id}\cap C_7$ is as shown in the next figure:
\vskip3mm
\begin{figure}[ht]
\begin{pspicture}(-3,-.5)(5,5)
      \psdots[dotsize=4pt](0,0)(1.071,0.214)
      \psline[linewidth=.01pt,linestyle=dashed,arrowsize=0.1]{->}(0,0)(4.16,5.2)
      \psline[linewidth=.01pt,linestyle=dashed,arrowsize=0.1]{->}(1.071,0.214)(4.7,4.75)
      \psline[linewidth=1.5pt]{->}(0,0)(1.32,3.96)
      \psline[linewidth=1.5pt]{->}(0,0)(2.475,1.65)
      \psline[linewidth=1.5pt]{->}(1.071,0.214)(4.843,1.629)
      \psline[linewidth=1.5pt](0,0)(1.071,0.214)
      \pspolygon[linewidth=.03pt](4,5)(1.2,3.6)(2.25,1.5)(4.5,1.5)(4.5,4.5)
      \uput{3pt}[215](0,0){$0$}
\end{pspicture}
\end{figure}
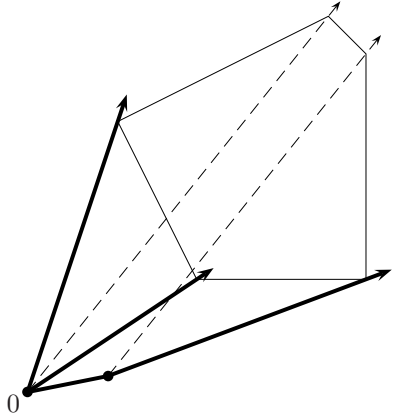
\centerline{Figure 2}
\vskip3mm
\noindent
{\bf 4. Intersections of $C_p$ and its faces with hyperplanes}
\vskip3mm
Certain classes of numerical semigroups correspond to the lattice points in the intersection of $C_p$ or some of its faces with hyperplanes.
\vskip2mm
\noindent
{\bf 4.1. Examples.}

\noindent
a) If $\mu(H)=(\mu_1,\dots,\mu_{p-1})$ for $H\in \mathfrak H_p$, then $g(H)=\sum_{i=1}^{p-1}\mu_i$ is the genus of $H$. Thus if $H_g$ is the hyperplane $\sum_{i=1}^{p-1}X_i=g$, then the lattice points in $C_p\cap H_g$ are in one-to-one correspondence with the semigroups in $\mathfrak H_p$ of genus $g$ and those of $C_p^0\cap H_g$ with the $H\in \mathfrak H_p$ of genus $g$ and maximal embedding dimension $p$.

\noindent
b) Given $q=i+np\ (i<p,i\ \text{and}\ p$ coprime) the lattice points in the intersection of $C_p$ with the hyperplane $X_i=n$ correspond bijectively to the $H\in \frak H_{pq}$ such that $q\in\text{ Ap}(H,p)$, and the lattice points in the intersection of $C_p$ with the half-space $X_i\le n$ to all $H\in \frak H_{pq}$. In fact, if $\mu(H)=(\mu_1,\dots,\mu_{p-1})$ and $\mu_i\le n$, then $i+\mu_ip\in H$ implies that also $q\in H$.
\vskip2mm
More generally, let $\alpha\in \mathbb N^{p-1}\setminus\{0\}$ be a primitive lattice point, that is, if $\alpha=(\alpha_1,\dots,\alpha_{p-1})$, then $\alpha_1,\dots,\alpha_{p-1}$ are relatively prime. For $x\in \mathbb R^{p-1}$ let $\alpha\cdot x$ denote the scalar product of $\alpha$ and $x$, and $H_n=H_n(\alpha):=\{x\in \mathbb R^{p-1}\vert \alpha\cdot x=n\}$ for $n\in \mathbb N$. We assume that no edge of $C_p^*$ is contained in the hyperplane $H_0$. This condition is satisfied in the Examples 4.1. If in 4.1b) there would be an edge vector
$\delta=(\mu_1,\dots,\mu_{i-1},0,\mu_{i+1},\dots,\mu_{p-1})$
there would be infinitely many lattice points on the edge of $C_p^*$
determined by $\delta$. These correspond to semigroups $H\in \frak H$
with $i\in H$. Since $i$ is prime to $p$ there exist only finitely
many such $H$, a contradiction.

Let $\bar S$ be the topological closure of a $d$-dimensional face $S$ of $C_p$ and let ${\it P}_{\alpha,n}$ resp. $\it P_{\alpha,n}^{\bar S}$ be the intersection of the hyperplane $H_n$ with the cone $C_p$ resp. $\bar S$. For $n\ge 1$ the sets $\it P_{\alpha,n}$ and ${\it P}_{\alpha,n}^{\bar S}$ are rational polytopes of dimension $p-2$ resp. $d-1$. We are interested in the numbers $f_{\alpha}(n)$ resp. $f_{\alpha}^{\bar S}(n)$ of lattice points in $P_{\alpha,n}$ and $P_{\alpha,n}^{\bar S}$ and in the number $f^0_{\alpha}(n)$ of lattice points in $C_p^0\cap P_{\alpha,n}$.
\vskip3mm\noindent
{\bf 4.2. Theorem.} Let ${\delta_1,\dots,\delta_t}$ be the canonical system of representatives for the edges of $C_p^*$ and ${\delta_1,\dots,\delta_r}$ the edge vectors contained in $\overline{S^*}$ where $S^*$ is the face of $C_p^*$ parallel to $S$. If $\bar S$ contains a lattice point $u\in \mathbb N^{p-1}$, then $f_{\alpha}^{\bar S}$ is a quasi-polynomial of degree $d-1$ with non-negative leading coefficient. The least common multiple of $\{\alpha\cdot\delta_i\}_{i=1,\dots,r}$ is a quasi-period of $f_{\alpha}^{\bar S}$. In particular $f_{\alpha}$ is a quasi-polynomial of degree $p-2$. Moreover for $n\ge\sum_{i=1}^{p-1}\alpha_i$
$$f^0_{\alpha}(n)=f_{\alpha}(n-\sum_{i=1}^{p-1}\alpha_i).$$
\vskip3mm

Remember that a function $f: \mathbb N\to \mathbb C, f\not\equiv 0$ with generating
function $H_f(T):=\sum_{n=0}^{\infty}f(n)T^n$ is called a {\it quasi-polynomial} of degree $d$ and quasi-period $N>0$ if $f$ is of the form
$$f(n)=c_d(n)n^d+c_{d-1}(n)n^{d-1}+\dots+c_0(n)\ (n\in \mathbb N)$$
where the $c_i:\mathbb N\to \mathbb C$ are periodic functions with integral period $N$ and $c_d(n)$ does not vanish identically. Equivalently, $f$ is a quasi-polynomial if there exist an integer $N>0$ and polynomials $f_0,\dots,f_{N-1}$ such that
$$f(n)=f_i(n)\ \text{if}\ n\equiv i\ \text{mod}\ N.$$
If there exists an integer $N>0$ and polynomials $P(T), Q(T)\in \mathbb C[T]\setminus \{0\}$ with $\text{deg}P(T)<\text{deg}Q(T)$ so that
$$H_f(T)=\frac{P(T)}{Q(T)}$$
and $\alpha^N=1$ for each zero $\alpha$ of $Q$, then $f$ is a quasi-polynomial with quasi-period $N$. Its degree is one less than the maximum pole order of the rational function $\frac{P(T)}{Q(T)}$ ([S2],4.4.1). $f\equiv 0$ is also considered as a quasi-polynomial.

\vskip3mm\noindent
{\bf Proof of Theorem 4.2.}
At first we show that $f_{\alpha}^{\bar S}$ is a quasi-polynomial. Since $H_n(\alpha)$ does not contain an edge of $\bar S^*$ the denominator $\prod_{i=1}^r(1-X^{\delta_i})$ of the generating function $H_{\bar S}$ (see Proposition 2.2)
does not vanish if we replace the $X_j$ by $T^{\alpha_j}$ with a variable $T$. Then
$$H_{\bar S}(T^{\alpha_1},\dots,T^{\alpha_{p-1}})=\frac{F_{\bar S}(T^{\alpha_1},\dots,T^{\alpha_{p-1}})}{\prod_{i=1}^r(1-T^{\alpha\cdot\delta_i})}=\sum_{n=0}^{\infty}f_{\alpha}^{\bar S}(n)T^n.$$
Set $P(T):=F_{\bar S}(T^{\alpha_1},\dots,T^{\alpha_{p-1}})$ and $Q(T):=\prod_{i=1}^r(1-T^{\alpha\cdot \delta_i})$. Since
deg$_{X_j}H_{\bar S}\le -1$ we have deg$P<\ $deg$Q$. Hence  $f_{\alpha}^{\bar S}(n)$ is a quasi-polynomial with quasi-period as stated. Its leading coefficient is non-negative since $f_{\alpha}^{\bar S}(n)\in\mathbb N$ for all $n\in \mathbb N$.

Now we show that $f_{\alpha}^{\bar S}$ has degree $d-1$. By assumption $\bar S$ contains a lattice point $u\in \mathbb N^{p-1}$. On the line through $u$ and the vertex $v$ of $C_p$ there is a lattice point $w$ not contained in $\bar S$, hence $-w\in \mathbb N^{p-1}$ and
$$u+\overline{S^*}\subset\bar S\subset w+\overline{S^*}.\leqno(1)$$
Let $\it P_{\alpha,n}^*:=\bar S^*\cap H_n$ for $n\in \mathbb N_+$. Then $\it P_{\alpha,1}^*$ is a convex polytope of dimension $d-1$ and $\it P_{\alpha,n}^*=n\cdot\it P_{\alpha,1}^*$. To this situation a theorem of Ehrhart ([S2],4.6.8) can be applied. It states that if $i(\it P_{\alpha,n}^*)$ is the number of lattice points in $\it P_{\alpha,n}^*$, then this function of $n$ is a quasi-polynomial of degree $d-1$.

In order to show this also for $f_{\alpha}^{\bar S}$ observe that (1) implies that for
$n\ge l:=\alpha\cdot u$ and $k:=\alpha\cdot w$ we have
$$u+\it P_{\alpha,n-l}^*\subset \bar S\cap H_n\subset w+\it P_{\alpha,n-k}^*.$$
Thus for large $n$ the quasi-polynomial $f_{\alpha}^{\bar S}$ is trapped by two quasi-polynomials of degree $d-1$, hence it has also the degree $d-1$.

The formula for $f_{\alpha}^0$ follows from Remark 2.1 after substituting $X_i=T^{\alpha_i}\ (i=1,\dots,p-1)$ and expanding into power series in $T$.\ $\Box$
\vskip3mm

In the situation of Example 4.1a) let $G(p,g)$ be the number of $H\in \mathfrak H_p$ with genus $g$. Then with $\alpha=(1,\dots,1)$ we have $G(p,g)=f_{\alpha}(n)$, and Theorem 4.2 tells us that $G(p,g)$, as a function of $g$, is a quasi-polynomial of degree $p-2$ with non-negative leading term. Moreover the least common multiple of $\{\alpha\cdot\delta_i\}_{i=1,\dots,r}$ is a quasi-period of $G(p,g)$.

$G^0(p.g):=f_{\alpha}^0(g)$ is the number of $H\in \mathfrak H_p$ having maximal embedding dimension $p$ and genus $g$. For $g\ge p-1$ we have
$$G^0(p,g)=G(p,g-(p-1)).\leqno(2)$$
If $G_{sym}(p,g)$ is the the number of symmetric $H\in \mathfrak H_p$ with genus $g$, then Proposition 3.1 and Theorem 4.2 imply that $G_{sym}(p,g)$ is a quasi-polynomial of degree
$\lfloor\frac{p}{2}\rfloor-1$.

In the situation of Example 4.1b) we can apply Theorem 4.2 to the functions $f_{e_i}^{\bar S}$ where $e_i$ is the i-th unit vector and $\bar S$ the closure of a face of $C_p$ containing a lattice point. Let $g_{e_i}(n):=\sum_{j=0}^nf_{e_i}(j)$ and $g_{e_i}^0(n):=\sum_{j=0}^nf^0_{e_i}(j)$, that is, the number of $H\in \mathfrak H_p$ (of maximal embedding dimension $p$) containing also $i+np$. In order to apply Theorem 4.2 also to $g_{e_i}$ and $g^0_{e_i}$ we need the following facts.

For functions $f: \mathbb N\to \mathbb C$ we consider the operators
$$E: f(n)\to f(n+1)\ \text{(shift})$$
$$\Delta: f(n)\to f(n+1)-f(n)\ \text{(difference})$$
$$\Sigma: f(n)\to \sum_{i=0}^nf(i)\ (\text{sum}).$$
It is easy to see that $f$ is a quasi-polynomial if and only if this is the case for $Ef, \Delta f$ or $\Sigma f$. Moreover deg$Ef$=deg$f$.
\vskip3mm
\noindent
{\bf 4.4. Lemma.} Let $f: \mathbb N\to \mathbb R, f\not\equiv 0$ be a quasi-polynomial of degree $d$.

\noindent
a) If $f$ is increasing, then $c_d(n)$ is constant.

\noindent
b) If $f\ge 0$, then deg$\Sigma f=d+1$.
\vskip3mm
\noindent
{\bf Proof.} a) Let $N$ be a quasi-period of $f=\sum_{k=0}^dc_k(n)n^k$ and $f_0,\dots,f_{N-1}$ the polynomials with $f(n)=f_i(n)$ for $i=0,\dots,N-1$. Then $c_d(i)$ is the coefficient of $t^d$ in $f_i(t)$ for $i=0,\dots,N-1$. With $k\in \mathbb N_+$ and $i\in \{0,\dots,N-1\}$ we have
$$f_0(kN)=f(kN)\le f(i+kN)=f_i(i+kN)\le f(N+kN)=f_0((k+1)N).$$
If follows that
$$c_d(0)=\lim_{k\to\infty}\frac{f_0(kN)}{(kN)^d}\le c_d(i)=\lim_{k\to\infty}\frac{f_i(i+kN)}{(kN)^d}\le\lim_{k\to\infty}\frac{f_0((k+1)N)}{(kN)^d}$$
$$=\lim_{k\to\infty}c_d(0)(\frac{k+1}{k})^d=c_d(0),$$
hence $c_d(0)=c_d(i)$ and $c_d$ is constant.

\noindent
b) If $c_d(n)=c_d$ is constant, then for $n\equiv i\ \text{mod}\ d$
$$f(n+1)-f(n)=c_d(n+1)^d-c_dn^d+g$$
with a quasi-polynomial $g$ of degree $\le d-1$. Hence deg$(\Delta f)\le d-1$ for an increasing $f$. As $\Sigma f$ is increasing for $f\ge 0$ and $\Delta\Sigma f=Ef$ we have
$$\text{deg}\Sigma f\ge\text{deg}\Delta(\Sigma f)+1=\text{deg}Ef+1=d+1.$$
Since $f\ge 0$ there exists $k>0$ so that for $n>k$ all polynomial functions $f_j(n)$ are increasing. Moreover there exists $i\in\{0,\dots,N-1\}$ such that $f_i(n)\ge f_j(n)$ for all $j$ and $n>k$, if $k$ is sufficiently large. Then for $n>k$
$$f(0)+f(1)+\dots+f(n)\le f(0)+\dots+f(k)+(n-k)f_i(n)=:h(n)$$
where $h(n)$ is a polynomial of degree $\le d+1$. Consequently deg$\Sigma f\le d+1$, and b) follows. $\Box$
\vskip3mm

Theorem 4.2 and Lemma 4.4 imply
\vskip3mm\noindent
{\bf 4.5. Proposition.} Let $i\in\{1,\dots,p-1\}$ be prime to $p$ and $N(p,i+np):=g_{e_i}(n)$ resp. Medim$(p,i+np):=g^0_{e_i}(n)$ the number of $H\in \mathfrak H_p$ (of maximal embedding dimension $p$) with $i+np\in H$. Then the functions $N(p,i+np)$ and $\text{Medim}(p,i+np)$ of the variable $n$ can
be expressed as quasi-polynomials of degree $p-1$ having the same
highest coefficient, which is independent of $n$. More precisely for each $q>p$ which is prime to $p$
$$\text{Medim}(p,q)=N(p,q-p).$$

For the last formula note that $f^0_{e_i}(0)=0, f^0_{e_i}(n)=f_{e_i}(n-1)$ for $n\ge 1$ by 4.2, hence Medim$(p,i+np)=g^0_{e_i}(n)=\sum_{k=0}^nf^0_{e_i}(k)=\sum_{k=0}^{n-1}f_{e_i}(k)=
g_{e_i}(n-1)=N(p,i+(n-1)p)$.

Similarly as in Proposition 4.5, if $H(p,g)$ denotes the number of $H\in \mathfrak H_p$ of genus $\le g$, then $H(p,g)$ is a quasi-polynomial of $g$ with degree $p-1$ and constant highest coefficient.
\vskip3mm \noindent

As in Proposition 3.1 let $\bar S_1,\dots,\bar S_t$ be the closures of the faces of dimension $\lfloor\frac{p}{2}\rfloor$ of $C_p$ whose lattice points correspond to the symmetric $H\in \mathfrak H_p$. Let Sym$(p,i+np)$ be the number of symmetric $H\in \mathfrak H_p$ containing $i+np$. Application of Theorem 4.2 and Lemma 4.4 to $\Sigma f_{e_i}^{\bar S_j}\ (j=1,\dots,t)$ yields
\vskip3mm\noindent
{\bf 4.6. Proposition.} $\text{Sym}(p,i+np)$ is, as a function of $n$, a quasi-polynomial of degree $\lfloor\frac{p}{2}\rfloor$ whose leading coefficient is independent of $n$.
\vskip3mm\noindent
{\bf 5. Asymptotic estimates and recursion formulas}
\vskip3mm

For an integer $q$ which is prime to $p$ let $\mathfrak H_{pq}$ be the set of all $H\in \mathfrak H_p$ with $q\in H$. We are interested in the functions $N(p,q), \text{Sym}(p,q)$ and Psym$(p,q)$ where $N(p,q)$ is the number of elements of $\mathfrak H_{pq}$,  Sym$(p,q)$ resp. Psym$(p,q)$ the number of symmetric resp. pseudo-symmetric $H\in \mathfrak H_{pq}$.

As in [KKW] and [KW] we associate with each $H\in \frak H_{pq}\ (q>p)$ a certain lattice path in the plane which is contained in the triangle $\Delta_0$ bounded by the line $g_0: p(X+1)+q(Y+1)=pq$ and the coordinate axes, starts on the $Y$-axis, ends on the $X$-axis and has only right or downward steps. In the following the word "lattice path" always means such a path. The set of all lattice paths for given $p,q$ is called the $(p,q)$-system.

The lattice path belonging to $H\in \frak H_{pq}$ is constructed as follows.
A semigroup $H\in \frak H_{pq}$ is obtained from $H_{pq}=<p,q>$ by closing some of its gaps. Each such gap $\gamma$ can be written
$$\gamma=pq-(a+1)p-(b+1)q\leqno (1)$$
with a unique $(a,b)\in \Delta_0$. The set $L_H$ of all these
$(a,b)$ is bounded by a lattice path and the coordinate axes. It is by definition
the path associated to $H$.

\vskip3mm
\pspicture*(-1.1,-1.1)(9,6.1)
\psgrid[gridlabels=0,subgriddiv=2,gridwidth=0.001pt,subgridwidth=0.001pt,gridcolor=lightgray](-1.1,-1)(9,6.1)
\psline[linewidth=1.5pt]{cc-cc}(-.5,5.5)(8,-.5)
\psline[linewidth=.5pt](0,0)(0,7)
\psline[linewidth=.5pt](0,0)(9,0)
\psline[linewidth=0.9pt](0,3.5)(0,2.5)
\psline[linewidth=0.9pt](0,2.5)(1,2.5)
\psline[linewidth=0.9pt](1,2.5)(1,2)
\psline[linewidth=0.9pt](1,2)(1.5,2)
\psline[linewidth=0.9pt](1.5,2)(1.5,1.5)
\psline[linewidth=0.9pt](1.5,1.5)(3,1.5)
\psline[linewidth=0.9pt](3,1.5)(3,1)
\psline[linewidth=0.9pt](3,1)(4,1)
\psline[linewidth=0.9pt](4,1)(4,0) \psdots*[dotsize=3pt](-.5,5.5)
\psdots*[dotsize=3pt](0,3.5) \psdots*[dotsize=3pt](1,2.5)
\psdots*[dotsize=3pt](1.5,2) \psdots*[dotsize=3pt](3,1.5)
\psdots*[dotsize=3pt](4,1) \psdots*[dotsize=3pt](4,0)
\uput{2pt}[45](1,1){$L$}\uput{3pt}[45](2,2){$\Delta_0$}
\psdots*[dotsize=3pt](8,-.5)
\uput{2pt}[45](2.5,3.6){$g_0$}
\uput{4pt}[180](0,3.5){$P_0$} \uput{4pt}[70](1,2.5){$P_1$}
\uput{4pt}[45](4,1){$P_{m-1}$} \uput{4pt}[270](4,0){$P_m$}
\endpspicture

\noindent
Figure 3

We use the notation $(P_0,P_1,\dots,P_m)$ for lattice paths where $P_0$ is the point where it starts, $P_m$ the point where it ends, and the other $P_i$ are the points where after a right step a downward step follows.
As is seen from
(1) a downward step in the lattice path of $H$ means for the
corresponding $\gamma$ an addition of q, a right step a subtraction
of $p$. If $(a,b)$ is a point of the lattice
path, then the corresponding $\gamma$ is the smallest element of $H$
in the residue class of $-(b+1)q$ modulo $p$, hence an element of
Ap($H,p)$, and different $b$ give different elements of the Ap\'ery
set. Thus the lattice path is given by the Ap\'ery set and
conversely also determines this set. This also indicates the
relation to the points of the polyhedral cone $C_p$.
\noindent
To $H=H_{pq}$ we may associate the empty path and the empty set $L_H$.

Not every lattice path as above belongs to a semigroup. For an arbitrary lattice path let $L$
be the set of lattice points in the area bounded by the path and the coordinate axes.
In order that $L=L_H$ with an $H\in \frak H_{pq}$ the following conditions must be satisfied ([KKW])

\noindent
a) For $(a,b),(a',b')\in L$ with $a+a'\ge q-1$ also $(a+a'-q+1,b+b'+1)$ must be in $L$.

\noindent
b) For $(a,b),(a',b')\in L$ with $b+b'\ge p-1$ also $(a+a'+1,b+b'-p+1)$ must be in $L$.

Lattice paths whose corresponding set $L$ satisfied a) and b) were
called {\it admissible}, and their number was denoted by $L(p,q)$.
\vskip3mm
\noindent
{\bf 5.1. Lemma.} $L(p,q)$ as a function of $p$ and of $q$ is increasing.
\vskip3mm
\noindent
{\bf Proof.} Let $(\tilde p,\tilde q)$ be another pair of relatively prime integers with $\tilde p\ge p,\tilde q\ge q$. Let $\tilde\Delta_0$ be the triangle bounded by the line $\tilde p(X+1)+\tilde q(Y+1)=\tilde p\tilde q$ and the axes. Clearly $\Delta_0\subset \tilde\Delta_0$ so that any lattice path in $\Delta_0$ is also one in $\tilde\Delta_0$. Let $L\subset\Delta_0$ be the set of lattice points corresponding to it. If it satisfies the admissibility conditions a) and b) above, then they are also satisfied in the $(\tilde p,\tilde q)$-system: Let $(a,b),(a',b')\in L$ and $a+a'\ge q-1$. If $(a+a'-q+1,b+b'+1)\in L$, then also $(a+a'-\tilde q+1,b+b'+1)=(a+a'-q+1,b+b'+1)-(\tilde q-q,0)\in L$. If $b+b'\ge p-1$ the proof is analogous.\ $\Box$
\vskip3mm

If $q=i+np$ and $g_{e_i}$ is the quasi-polynomial studied in Section 4 we have
$$N(p,q)=L(p,q)+1=g_{e_i}(n)\leqno (2)$$
where the 1 comes from $H_{pq}$ or the empty lattice path. Proposition 4.5 and Lemma 5.1 imply
\vskip3mm\noindent
{\bf 5.2. Proposition.} $\lim_{q\to\infty}\frac{N(p,q)}{q^{p-1}}$ and $\lim_{q\to \infty}\frac{\text{Medim}(p,q)}{q^{p-1}}$
exist and are equal.
\vskip3mm

The above limits give asymptotic estimates of how many $H\in \mathfrak H_{pq}$ (of maximal embedding dimension $p$) exist. With a somewhat different approach it is shown in [HW] that also for an arbitrary $q$ which is prime to $p$ the function $N(p,q)$ is a quasi-polynomial in $q$ of degree $p-1$ and estimates of its (constant) highest coefficients are given, i.e. of the above limits.

It was shown in [KW] that the admissible lattice paths starting at the point $(0,p-2)$
are in one-to-one correspondence with the semigroups of $\frak H_{p,q-p}$ (if $q-p<p$ exchange $p$ and $q-p$).
Thus we have the recursion formula
$$N(p,q)=N_{pq}+N(p,q-p)+1\qquad (q>p)\leqno (3)$$
where $N_{pq}$ is the number of admissible lattice paths starting at
$(0,b)$ with $b\in \mathbb N, b\le p-3$.

Also for the numbers $\text{Sym}(p,q)$ resp. $\text{Psym}(p,q)$ of symmetric (pseudo-symmetric) $H\in \frak H_{pq}$ one has recursion formulas
$$\text{Sym}(p,q)=S_{pq}+\text{Sym}(p,q-p)+1\qquad(q>p)\leqno(4)$$
$$\text{Psym}(p,q)=P_{pq}+\text{Psym}(p,q-p)\qquad (q>p)\leqno(5)$$
where $S_{pq}$ (resp. $P_{pq}$) is the number of admissible lattice
paths starting at a point $(0,b)$ with $b\le p-3$ and defining a
symmetric (pseudo-symmetric) semigroup. We shall use the formulas
(3)-(5) in the next section to derive explicit recursion formulas for $N(p,q),
\text{Sym}(p,q)$ and $\text{Psym}(p,q)$ in the
cases $p=3, p=4$.
\hfill\eject\noindent
{\bf 6. Examples}
\vskip3mm
The function $G(p,g)$ tells us in how
many ways we can remove $g$ numbers from $\mathbb N$, where $0$ and
$p$ are not removed, so that the remaining set is additively closed.
In [Ku], Appendix C various generating functions for $\alpha=(1,\dots,1)$ and $p\le 5$ are listed. The following explicit formulas can be easily derived:
$$G(3,g)=\lfloor\frac{g}{3}\rfloor+1$$ $$G(4,g)=\frac{1}{12}g^2+\frac{1}{2}g+\begin{cases}1\quad \text{if}\ g\equiv 0\ \text{mod}\ 6\cr
\frac{5}{12}\ \text{if}\ g\equiv 1\ \text{mod}\ 6\cr
\frac{2}{3}\quad \text{if}\ g\equiv 2\ \text{mod}\ 6\cr
\frac{3}{4}\quad \text{if}\ g\equiv 3\ \text{mod}\ 6\cr
\frac{2}{3}\quad \text{if}\ g\equiv 4\ \text{mod}\ 6\cr
\frac{5}{12}\ \text{if}\ g\equiv 5\ \text{mod}\ 6\end{cases}=
\lfloor\frac{1}{12}g^2+\frac{1}{2}g\rfloor+1$$

$$G_{sym}(3,g)=\begin{cases}1\quad \text{if}\ g\not\equiv 2\ \text{mod}\ 3\cr
0\quad \text{if}\ g\equiv 2\ \text{mod}\ 3\end{cases}\  G_{sym}(4,g)=\lfloor\frac{g}{3}\rfloor+1.$$
For $p=3$ see also Fig.1. For $G^0(p,g)$ see formula (2) of Section 4.
\vskip3mm
$$\begin{matrix} g & G(3,g)&G^0(3,g)&G_{sym}(3,g)&G(4,g)&G^0(4,g)&G_{sym}(4,g)\cr
\noalign{\hrule}
0&1&0&1&1&0&1\cr
1&1&0&1&1&0&1\cr
2&1&1&0&2&0&1\cr
3&2&1&1&3&1&2\cr
4&2&1&1&4&1&2\cr
5&2&2&0&5&2&2\cr
6&3&2&1&7&3&3\cr
7&3&2&1&8&4&3\cr
8&3&3&0&10&5&3\cr
\end{matrix}$$
\vskip3mm
The following formulas were communicated to us by H. Knebl:
$$G(5,g)=\frac{1}{135}g^3+\frac{4}{45}g^2+R(i)$$
$$G_{sym}(5,g)=\begin{cases}\frac{1}{6}g+S(i)\quad\ g\not\equiv 3\ \text{mod}\ 5\cr
\qquad 0\qquad\ \ g\equiv 3\ \text{mod}\ 5\end{cases}$$
with $R(i)$ and $S(i)$, depending on $i\equiv g\ \text{mod}\ 30$, as in the table below.

$$\begin{matrix} i & R(i)& S(i)\cr
\noalign{\hrule}
0&(7/15)g+1&1\cr
1&(1/3)g+77/135&5/6\cr
2&(19/45)g+20/27&2/3\cr
3&(3/10)g+1/10&-\cr
4&(2/5)g+68/135&4/3\cr
5&(29/90)g+13/54&1/6\cr
6&(13/30)g+3/5&1\cr
7&(11/30)g+29/54&5/6\cr
8&(16/45)g+91/135&-\cr
9&(3/10)g+7/10&1/2\cr
10&(7/15)g+28/27&4/3\cr
11&(13/45)g+28/135&1/6\cr
12&(7/15)g+4/5&1\cr
13&(3/10)g-53/270&-\cr
14&(16/45)g+37/135&2/3\cr
15&(11/30)g+1/2&1/2\cr
16&(13/30)g+131/135&4/3\cr
17&(29/90)g+119/270&1/6\cr
18&(2/5)g+4/5&-\cr
19&(3/10)g+109/270&5/6\cr
20&(19/45)g+20/27&2/3\cr
21&(1/3)g+1/5&1/2\cr
22&(7/15)g+113/135&4/3\cr
23&(23/90)g-7/270&-\cr
24&(2/5)g+4/5&1\cr
25&(11/30)g+29/54&5/6\cr
26&(7/18)g+82/135&2/3\cr
27&(11/30)g+1/2&1/2\cr
28&(2/5)g+68/135&-\cr
29&(23/90)g+47/270&1/6\end{matrix}$$
In [Ka], table 1 an extensive list with values of the function $S(m,g)$ counting the semigroups with multiplicity $m$ and genus $g$ is given. If $m$ is a prime number, then $S(m,g)=G(m,g)$ for large $g$, since there are only finitely many $H\in \mathfrak H_m$ of multiplicity $<m$.

In the following the recursion formulas (3)-(5) from Section 5 will be used.

\noindent
I) For $\bf{p=3}$, in order to determine $N_{3q}, S_{3q}$ and $P_{3q}$, we have only to consider admissible lattice paths on the $X$-axis. These are the paths ending at $(j,0)$ with $j\le \lfloor\frac{q}{2}\rfloor-1$. They correspond  to the semigroups $H_j=<3,q,2q-3(j+1)>$ where $H_j$ is symmetric if and only if $q$ is even and $j=\frac{q}{2}-1$ and $H_j$ is pseudo-symmetric if and only if $j=0$ or $q$ is odd and $j=\frac{q-1}{2}-1$ ([KW], Example 2.5b). Thus $N_{3q}=\lfloor\frac{q}{2}\rfloor$ and
$$N(3,q)=N(3,q-3)+\lfloor\frac{q}{2}\rfloor+1\qquad (q>3),\leqno(1)$$
and for $q\ge 4$
$$\text{Sym}(3,q)=\text{Sym}(3,q-3)+\begin{cases}2\quad \text{if $q$ is even}\cr 1\quad\text{if $q$ is odd}\end{cases}\leqno(2)$$

$$\text{Psym}(3,q)=\text{Psym}(3,q-3)+\begin{cases}1\quad\text{if $q$ is even}\cr 2\quad\text{if $q$ is odd}\end{cases}.\leqno(3)$$
Clearly $\text{Sym}(3,1)=1$ and $\text{Sym}(3,2)=2$,
further $\text{Psym}(3,1)=\text{Psym}(3,2)=0$.
For Medim$(p,q)$, thanks to Proposition 4.5, no further discussion is necessary.

$$\begin{matrix} q & N(3,q)& \text{Medim}(3,q)&\text{Sym}(3,q)&\text{Psym}(3,q)\cr
\noalign{\hrule}
1&1&0&1&0\cr
2&2&0&2&0\cr
4&4&1&3&1\cr
5&5&2&3&2\cr
7&8&4&4&3\cr
8&10&5&5&3\cr
10&14&8&6&4\cr
11&16&10&6&5\cr
13&21&14&7&6\cr
14&24&16&8&6\end{matrix}$$
The numbers in the columns can also be found by using figure 1.
By what was said there we have
$N(3,q)=\text{Medim}(3,q)+\text{Sym}(3,q)$, and by Proposition 4.5 Medim$(3,q)=N(3,q-3)$.
\vskip2mm
\noindent
{\bf 6.1. Proposition.} a) $\lim_{q\to\infty}\frac{N(3,q)}{q^2}=\lim_{q\to \infty}\frac{\text{Medim}(3,q)}{q^2}=\frac{1}{12}$.

\noindent
b) $\lim_{q\to\infty}\frac{\text{Sym}(3,q)}{q}=\lim_{q\to\infty}\frac{\text{Psym}(3,q)}{q}=\frac{1}{2}$.
\vskip3mm\noindent
{\bf Proof.} a) By [KW], Example 3.5 we have $N(3,q)=\lfloor\frac{q^2}{12}+\frac{q}{2}\rfloor+1$ which implies a).
Alternately, from the recursion formula (1) follows for $q\ge 7$ that
$$N(3,q)-N(3,q-6)=\lfloor\frac{q}{2}\rfloor+\lfloor\frac{q-3}{2}\rfloor+2=q$$
$$N(3,q-6(i-1))-N(3,q-6i)=q-6(i-1)\ (i\ge 1)$$
hence
$$N(3,q)-N(3,q-6\lfloor\frac{q}{6}\rfloor)=\lfloor\frac{q}{6}\rfloor q-6\sum_{i=0}^{\lfloor\frac{q}{6}\rfloor-1}i=$$
$$\frac{q^2}{12}+\frac{q}{2}-
3(\frac{q}{6}-\lfloor\frac{q}{6}\rfloor+1)(\frac{q}{6}-\lfloor\frac{q}{6}\rfloor)$$
which also gives a).

\noindent
b) By (2) and (3) for $q\ge 7$
$$\text{Sym}(3,q)-\text{Sym}(3,q-6)=3=\text{Psym}(3,q)-\text{Psym}(3,q-6),$$
hence
$$\text{Sym}(3,q)-\text{Sym}(3,q-6\lfloor\frac{q}{6}\rfloor)=3\lfloor\frac{q}{6}\rfloor=
\text{Psym}(3,q)-\text{Psym}(3,q-6\lfloor\frac{q}{6}\rfloor),$$
from which b) follows.\ $\Box$
\vskip2mm\noindent
II) For $\bf{p=4}$, in order to determine $N_{4q}, S_{4q}$ and $P_{4q}$, we have to investigate the lattice paths starting at $(0,1)$ or $(0,0)$ and the corresponding semigroups.

By [KW], 2.6 there are ${{2+q'}\choose{2}}-1=\frac{1}{8}(q^2+4q+3)-1\ (q'=\frac{q-1}{2})$ lattice paths in the rectangle $\bf R$ with the corners $(0,0),(0,1),(q'-1,1),(q'-1,0)$
and all are admissible. To find $N_{4q}$ we have also to determine the admissible paths starting at $P_0=(0,1)$ and ending
at a point $(q'+i,0)\ (i=0,\dots,q-2-\lfloor\frac{q}{4}\rfloor-q')$.

Such a lattice path is admissible if and only if it contains $(2i,1)$ and avoids $(q'-i,1)$. Thus it is of the form $(P_0,P_1,P_2)$ with $P_1=(j,1),P_2=(q'+i,0)\ (2i\le j<q'-i)$ and defines the semigroup
$$G_{ij}:=<4,q,2q-4(j+1),3q-4(q'+i+1)>.$$
Given $i$ there are $q'-3i$ such paths, the number of possible downward steps, where $i\le\lfloor\frac{q'}{3}\rfloor=\lfloor\frac{q-1}{6}\rfloor$. Note that the last condition implies that
$$i\le q-2-\lfloor\frac{q}{4}\rfloor-q'\leqno(4)$$
so that we stay below the line $g_0$.
Altogether there are
$$\sum_{i=0}^{\lfloor\frac{q-1}{6}\rfloor}(q'-3i)=
(\lfloor\frac{q-1}{6}\rfloor+1)\frac{q-1}{2}-\frac{3}{2}(\lfloor\frac{q-1}{6}\rfloor+1)\lfloor\frac{q-1}{6}\rfloor$$
such lattice paths (semigroups $G_{ij}$), and therefore for $q>4$
$$N(4,q)=
\frac{1}{8}(q^2+4q+3)+(\lfloor\frac{q-1}{6}\rfloor+1)(\frac{q-1}{2}-\frac{3}{2}\lfloor\frac{q-1}{6}\rfloor)
+N(4,q-4).\leqno(5)$$
Clearly $N(4,1)=1, N(4,3)=4$.

In order to determine the numbers $S_{4q}$ and $P_{4q}$ we first check which of the semigroups $G_{ij}$ are symmetric resp. pseudo-symmetric. This is done by considering their genus and Frobenius number.

$G_{ij}$ is obtained from $H_{4q}$ by closing $q'+i+j+2$ of its gaps, hence $G_{ij}$ has genus
$$g(G_{ij})=\frac{3}{2}(q-1)-(q'+i+j+2)=2q'-(i+j+2).$$
The parallel $4(X+1)+q(Y+1)=3q+4$ to $g_0$ through $(0,2)$ passes also through $(\frac{q}{4},1)$ and cuts the $X$-axis at $(\frac{q}{2},0)$. If $j\ge \lfloor\frac{q}{4}\rfloor$, then by [KW], 3.1 the point $(0,2)$ corresponds to the Frobenius number of $G_{ij}$, hence
$$F(G_{ij})=3(q-1)-1-2q=q-4=2q'-3$$
and
$$2g(G_{ij})-F(G_{ij})-1=2(q'-i-j-1),$$
an even number. Therefore $G_{ij}$ is not pseudo-symmetric, and symmetric if and only if $j=q'-i-1$. Furthermore $q'-i-1\ge\lfloor\frac{q}{4}\rfloor$ by (4). Thus $\lfloor\frac{q-1}{6}\rfloor+1$ of the $G_{ij}$ are symmetric for $q\not\equiv 1\ \text{mod}\ 6$, and $\lfloor\frac{q-1}{6}\rfloor$ are symmetric for $q\equiv 1\ \text{mod}\ 6$. In case $j<\lfloor\frac{q}{4}\rfloor$ the point $(j+1,1)$ corresponds to $F(G_{ij})$, hence
$$F(G_{ij})=3(q-1)-1-4(j+1)-q=4q'-4j-6$$
and
$$2g(G_{ij})-F(G_{ij})-1=2j-2i+1.$$
In this case $G_{ij}$ is not symmetric, and pseudo-symmetric only for $j=i$. Since $j\ge 2i$ this implies $i=j=0$.

According to [KW] 2.8, of the semigroups $\tilde G_{ij}$ corresponding
to lattice paths in the rectangle $\bf R$ exactly $q'$ are symmetric and just one is pseudo-symmetric. Altogether for $q\ge 7$
$$S_{4q}=\begin{cases}\lfloor\frac{q-1}{6}\rfloor+\frac{q-1}{2}\qquad \text{if}\ q\equiv 1\ \text{mod}\  6\cr
\lfloor\frac{q-1}{6}\rfloor+\frac{q-1}{2}+1\ \text{if}\ q\not\equiv 1\ \text{mod}\ 6 \end{cases}$$
and $P_{4q}=2.$ Thus
$$\text{Sym}(4,q)=\text{Sym}(4,q-4)+1+\begin{cases}\lfloor\frac{q-1}{6}\rfloor+
\frac{q-1}{2}\qquad \text{if}\  q\equiv 1\ \text{mod}\ 6\cr
\lfloor\frac{q-1}{6}\rfloor+\frac{q-1}{2}+1\ \text{if}\ q\not\equiv 1\ \text{mod}\ 6 \end{cases}$$
$$\text{Psym}(4,q)=\text{Psym}(4,q-4)+2.$$
It is easy to see that $\text{Sym}(4,3)=3, \text{Sym}(4,5)=5$,
$\text{Psym}(4,3)=1, \text{Psym}(4,5)=2$.

$$\begin{matrix} q & N(4,q)& \text{Medim}(4,q)&\text{Sym}(4,q)&\text{Psym}(4,q)\cr
\noalign{\hrule}
1&1&0&1&0\cr
3&4&0&3&1\cr
5&9&1&5&2\cr
7&17&4&8&3\cr
9&29&9&11&4\cr
11&45&17&15&5\cr
13&66&29&19&6\cr
15&93&45&25&7
\end{matrix}$$

\vskip3mm\noindent
{\bf 6.2. Proposition.} $\lim_{q\to\infty}\frac{N(4,q)}{q^3}=\lim_{q\to\infty}\frac{\text{Medim}(4,q)}{q^3}=\frac{1}{72}$.
\vskip3mm\noindent
{\bf Proof.} We know by 5.2 that the limits exists and are equal. Therefore it is enough to consider $q$ of the form $q=4m+1\ (m>0)$. By formula (5) we obtain for $q\ge 5$
$$N(4,q)-N(4,q-4)=\frac{1}{8}(q^2+4q+3)+(\lfloor\frac{q-1}{6}\rfloor+1)(\frac{q-1}{2}
-\frac{3}{2}\lfloor\frac{q-1}{6}\rfloor)$$
$$=\frac{8}{3}m^2+g(m)$$
with a quasi-polynomial $g$ of degree $1, g(0)=1$. Therefore with the sum operator $\Sigma$
$$N(4,q)=\frac{8}{3}\sum_{j=1}^m j^2+(\Sigma g)(m).$$
By Fibonacci's formula (Liber abaci 1202) $\sum_{j=1}^m j^2=\frac{1}{6}m(m+1)(2m+1)$ this implies
$$N(4,q)=\frac{8}{9}m^3+h(m)=\frac{1}{72}(q-1)^3+h(m)$$ with a quasi-polynomial $h$ of degree $2$, and the assertion about the limits follows.\ $\Box$
\vskip3mm

$N(4,q)$ and $N(5,q)$ have been explicitly computed by H. Knebl (see [HW], Example 4.3).
\vskip3mm
{\bf References}
\vskip3mm\noindent\footnotesize
{[A] R.Ap\'ery. Sur les branches superlin\'eaires des courbes algebriques. C.R.Acad.
Sci. Paris 222, 1198-1200 (1946)
\vskip3mm
\noindent
[BGP] V. Blanco, P. A. Garc\'ia-S\'anchez, J. Puerto. Counting numerical semigroups with short generating functions. Intern. Journ. of Algebra and Computation 21 (2011) 1-18
\vskip3mm\noindent
[HW] M. Hellus, R. Waldi. On the number of numerical semigroups containing two coprime integers $p$ and $q$.
\vskip3mm\noindent
[Ka] N. Kaplan, Counting numerical semigroups by genus and some cases of a question of Wilf. Journ. Pure and Appl. Algebra 216 (2012) 1016-1032
\vskip3mm \noindent
[KKW] H. Knebl, E. Kunz, R. Waldi. Weierstra\ss\ semigroups and nodal curves of type
p,q. Journal of Algebra 348 (2011) 315-335
\vskip3mm \noindent
[Ku] E. Kunz. \"Uber die Klassifikation numerischer Halbgruppen. Regensburger mathemati-sche Schriften 11. Regensburg 1987. Available at

\noindent
http://www-nw.uni-regensburg.de/kue22107.mathematik.uni-regensburg.de/index.htm
\vskip3mm\noindent
[KW] E. Kunz, R. Waldi. Geometrical illustration of numerical semigroups and of some of their invariants. To appear in Semigroup Forum. DOI:10.1007/s00233-014-9599-7 (June 2014)
\vskip3mm\noindent
[RG] J. C. Rosales, P. A. Garc\'ia-S\'anchez. Numerical Semigroups. Springer 2009
\vskip3mm\noindent
[RGGB] J. C. Rosales, P. A. Garc\'ia-S\'anchez, J. I. Garc\'ia-Garc\'ia, M. B. Branco. Systems of inequalities and numerical semigroups. J. London Math. Soc. 65 (2002) 611-623
\vskip3mm\noindent
[S1] R. P. Stanley. Combinatorics and Commutative Algebra. Progress in Math. 41. Birkh\"auser Boston-Basel-Stuttgart 1983
\vskip3mm\noindent
[S2] ---- Enumerative Combinatorics, Vol.1 (sec. edition). Cambridge Studies in Advanced Mathematics 49. Cambridge University Press 2012}

\end{document}